\documentclass{ecgd-l}

\usepackage{amssymb}

\usepackage{graphicx}

\newcommand{\R}{\mathbb{R}}
\newcommand{\La}{\Lambda}
\newcommand{\C}{\mathbb{C}}
\newcommand{\Z}{\mathbb{Z}}

\newtheorem{theorem}{Theorem}[section]

\theoremstyle{definition}
\newtheorem{definition}[theorem]{Definition}

\theoremstyle{remark}
\newtheorem{remark}[theorem]{Remark}

\numberwithin{equation}{section}

\begin{document}

\title{Translation covers of some triply periodic platonic surfaces}


\author{Jayadev S.~Athreya}
\address{University of Washington, Seattle, WA 98195}
\curraddr{}
\email{jathreya@uw.edu}
\thanks{}

\author{Dami Lee}
\address{University of Washington, Seattle, WA 98195}
\curraddr{}
\email{damilee@uw.edu}
\thanks{The second author was supported by NSF Grant No. DMS-1440140}


\date{}

\begin{abstract}
We study translation covers of several triply periodic polyhedral surfaces that are intrinsically Platonic. We describe their affine symmetry groups and compute the quadratic asymptotics for counting saddle connections and cylinders, including the count of cylinders weighted by area. The mathematical study of triply periodic surfaces was initiated by Novikov, motivated by the study of electron transport. The surfaces we consider are of particular interest as they admit several different explicit geometric and algebraic descriptions, as described, for example, in the second author's thesis.
\end{abstract}

\maketitle


\section{Introduction}
\label{sec: intro}
In this paper, we study translation covers of six intrinsically Platonic polyhedral surfaces. The objects of our interest are quotients of triply periodic polyhedral surfaces that arise in Lee~\cite{thesis}. We say a polyhedral surface $\Pi$ is triply periodic, if $\Lambda \Pi = \Pi$ for some rank-three lattice of translations $\Lambda \subset \R^3.$ We use $\Pi$ to denote the infinite surfaces, and $X:=\Pi / \Lambda$ for the smallest compact quotients. The objects of our interest are those whose underlying surfaces $X$ are compact with identifiable Riemann surface structures. 

The study of triply periodic surfaces arises in an important physical situation, known as Novikov's problem on understanding the geometry of sections of compact surfaces embedded in the $3$-torus $\mathbb T^3$. Novikov's problem is motivated by semiclassical motion of an electron on the (periodic) Fermi surface of a metal under a magnetic field~\cite{Zorich:novikov}. 

In this paper, we are interested in \emph{polyhedral} surfaces embedded in the $3$-torus $\mathbb T^3$, or, equivalently, triply ($\mathbb Z^3$)-periodic polyhedral surfaces in $\mathbb R^3$. To be clear, we are \emph{not} studying Novikov's problem, but rather, we are interested in the intrinsic flat geometry of these objects given by their polyhedral structure. The particular surfaces we study also have rich automorphism groups and admit explicit descriptions as algebraic curves, see Lee~\cite{thesis} for further details. 

We are interested in the counting and geometry of special geodesic trajectories, using ideas of Veech~\cite{veech} and Gutkin--Judge~\cite{gj} and others to compute explicit quadratic asymptotics of the growth rate of the number of such trajectories of bounded length. That these asymptotic constants (known as Siegel--Veech constants) can be computed is by now well-known, but we believe it is of interest to compute them explicitly in these cases, as these surfaces occur naturally in many different contexts.

\subsection{Background}
\subsubsection*{Polyhedral surfaces}
We denote a polyhedral surface by $\{p, q\}$ (using Schl{\"a}fli symbols) if it is tiled by regular Euclidean $p$-gons and all vertices are $q$-valent. A polyhedral surface is \textit{Platonic} if there are two types of symmetries (an order-$p$ symmetry about the center of a polygonal face and an order-$q$ symmetry about a vertex) as isometries of the ambient space so that group generated by these two symmetries acts transitively on the vertices, edges, and faces. 

\begin{figure}[htbp] 
	\centering
	\includegraphics[width=4.5in]{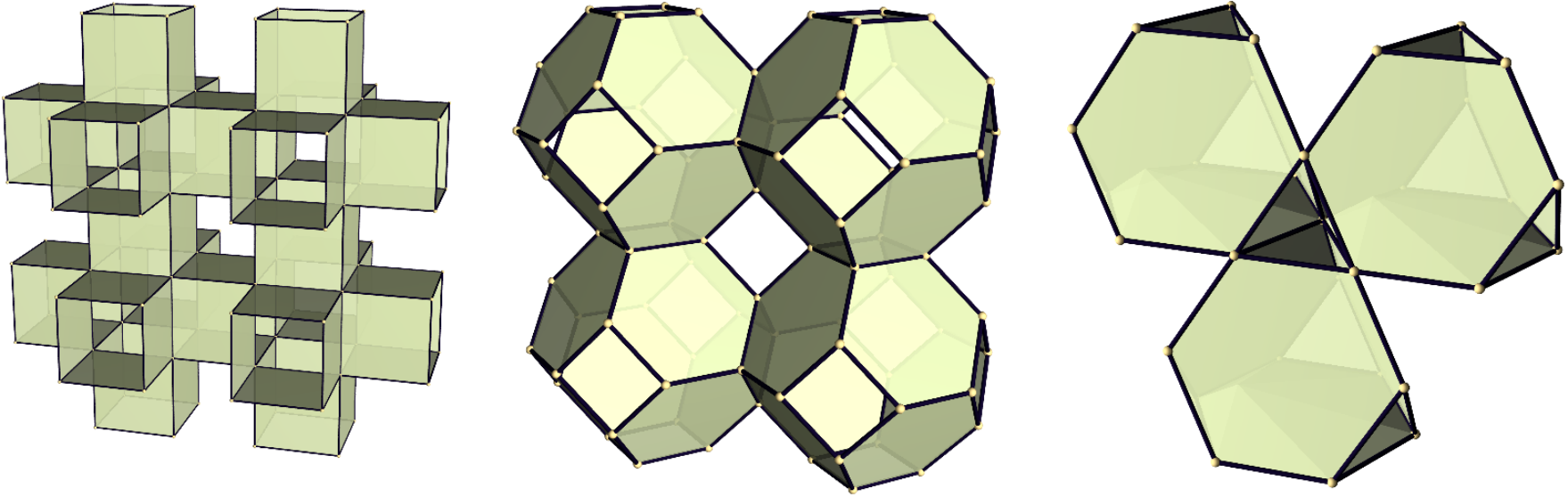}
	\caption{A subset of Mucube, Muoctahedron, and Mutetrahedron. Adapted from \cite{thesis}.}	
	\label{fig: platonic surfaces 1} 
\end{figure}

\begin{figure}[htbp] 
	\centering
	\includegraphics[width=4.5in]{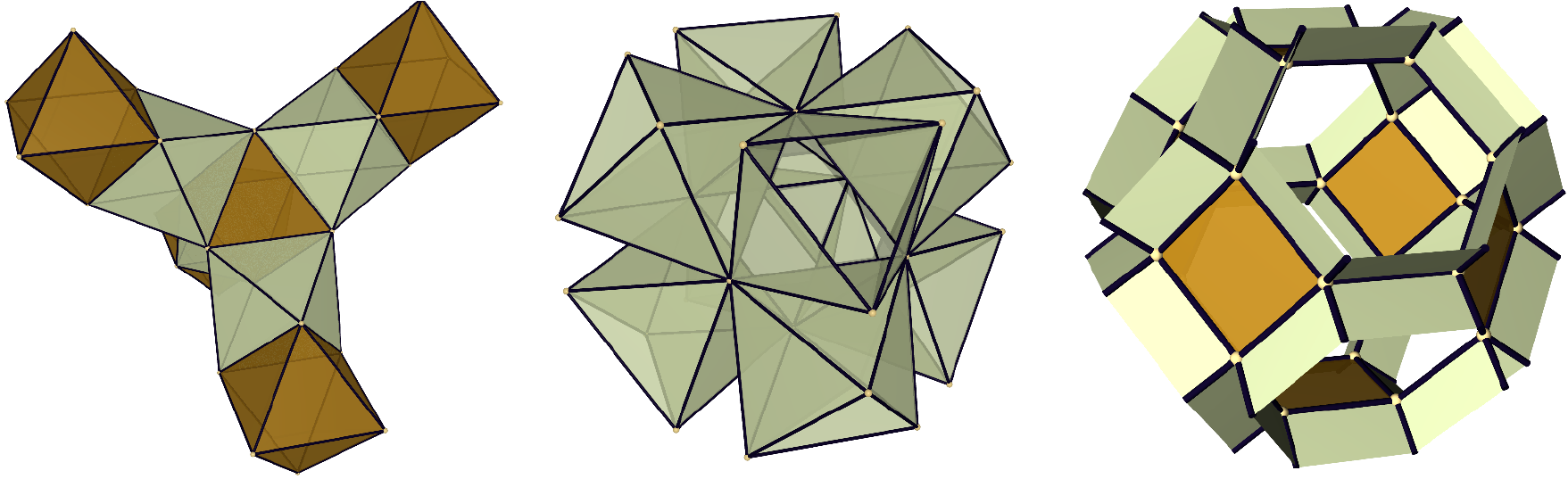}
	\caption{A subset of Octa-4, Octa-8, and Truncated Octa-8. Adapted from \cite{thesis}.}	
	\label{fig: rad} 
\end{figure}

We say that a surface is \emph{intrinsically Platonic} if the isometries of order-$p$ and $q$ are intrinsic and not necessarily Euclidean. That is, we disregard the embedding of the surface in Euclidean space and consider only the intrinsic metric on $X.$ Lee's classification~\cite{thesis} includes classical examples found by Coxeter and Petrie \cite{CP} such as Mucube $\{4, 6\},$ Muoctahedron $\{6,4\},$ and Mutetrahedron $\{6,6\}$ (Figure~\ref{fig: platonic surfaces 1}), and introduces examples that we call Octa-4 $\{3,8\},$ Octa-8 $\{3,12\},$ and Truncated Octa-8 $\{4,5\}$ (Figure~\ref{fig: rad}). The construction of these surfaces is described in Section~\ref{sec: platonic surfaces}. We remark that the genus of the underlying surface is three for Coxeter--Petrie's examples and Octa-4; and four for Octa-8 and Truncated Octa-8.

\subsubsection*{Translation surfaces}\label{sec:translation} A \emph{translation surface} is a collection of Euclidean polygons in $\R^2$ with parallel sides identified by translation. Equivalently, it is given by a pair $(X,\omega)$ where $X$ is a compact Riemann surface and $\omega$ is a holomorphic $1$-form (i.e., a section of the canonical bundle). In local coordinates, we can write $\omega = f(z)dz$. An order-$k$ zero of $\omega$ corresponds to a point with cone angle $2\pi(k+1)$. We say that $(X,\omega)$ lies in a stratum $\mathcal{H}_1(k_1, \ldots , k_n)$ where $k_i$ denotes the order of a zero of $\omega.$ Then $k_1 + \cdots + k_n = 2 g - 2$ where $g$ is the genus of $X.$ Integrating $\omega$ away from its zeros, we get an atlas of charts to $\mathbb C$ whose transition maps are translations. 

The group $GL^+(2,\R)$ acts on the moduli space of genus-$g$ translation surfaces. If the $GL^+(2,\R)$-orbit of the surface is closed, we call the $GL^+(2,\R)$-orbit a \emph{Teichm{\"u}ller curve.} The integral $\frac{i}{2} \int_X \omega \wedge \bar{\omega}$ represents the area of the translation surface, $\mbox{area}(\omega),$ and the group $SL(2, \R)$ acts on the moduli space of genus-$g$ unit area translation surfaces via $\R$-linear post-composition with charts. 

\subsubsection*{Saddle connections and cylinders}\label{sec:saddles}

A \textit{saddle connection} on a translation surface is a geodesic segment that connects two singular points (zeros of $\omega$) with no singular points in its interior. It is a \textit{closed saddle connection} if it connects a singular point to itself. We will be interested in counting saddle connections (organized by length) and understanding the existence of closed saddle connections. Given a saddle connection $\gamma$, its holonomy vector is given by $$z_{\gamma} = \int_{\gamma} \omega,$$ and we denote the set of all holonomy vectors by $\La_{\omega}$. This is a discrete subset of $\C,$ and it  varies equivariantly under the $SL(2, \mathbb R)$-action, so for $h \in SL(2, \R),$ $$h\La_{\omega}= \La_{h\omega}.$$ A closed geodesic $\eta$ not passing through a singular point is part of a \emph{cylinder}, and we denote the set of holonomy vectors of cylinders as $\La^{cyl}_{\omega}.$ Associated to each cylinder is the \emph{area} of the cylinder, $a(\eta).$

\subsubsection*{Translation covers}\label{sec:translationcover} Let $X$ be a Riemann surface and $k$ be a positive integer. A meromorphic $k$-differential on $X$ is a section of the $k$-th power of the canonical bundle $q = f(z) (d z)^k.$ In other words, integrating $\sqrt[k]{q}$ yields transition functions that are given by rotations by a multiple of $\frac{2\pi}{k}.$

Our polyhedral surfaces (as polygons in $\R^2$ with edge identifications) are not translation surfaces. Associated to a pair $(Y, \sigma)$, where $\sigma$ is a $k$-differential, is the \emph{translation cover} (also called \emph{spectral curve} or \emph{unfolding}), a translation surface $(X, \omega)$ where $X$ is a $k$-cover of $Y$ branched at the zeros of $\sigma$. Geometrically, it is constructed by taking $k$ copies of $(Y, \sigma),$ each rotated by $\frac{2\pi}{k}$ from the previous copy, where an edge $j$ of some copy of $Y$ is identified to an edge $j$ of some other copy by translation. A key point in our work is precisely identifying these covers in a variety of situations, and carefully computing the affine symmetry groups of the resulting translation surfaces, using key ideas from Gutkin--Judge~\cite{gj}. 

Moreover, we show that certain polyhedral surfaces and Platonic solids have common translation covers.

\begin{theorem}\label{thm:common covers}
\begin{enumerate}
\item The translation covers of the octahedron and the compact quotient of Octa-8 are isometric with their intrinsic metrics. 
\item The translation covers of the cube and compact quotient of Mucube are isometric with instrinsic metrics.
\end{enumerate}
\end{theorem}

It is natural to extend this question to more general cases.

Let $\Pi$ be a triply periodic polyhedral surface and $\Pi/\Lambda$ a compact Riemann surface equipped with a $k$-differential. For $k' | k,$ does there exist a triply periodic polyhedral surface $\Pi'$ such that $\Pi'/\Lambda'$ is a compact Riemann surface that is a $k'$-cover of $\Pi / \Lambda?$ If not, what are the constraints? Specifically, can one find a triply periodic polyhedral surface whose compact quotient is a 2-cover (or 5-, 10-) of the regular icosahedron?

\subsubsection*{Veech groups and counting results}\label{sec:lattice}
The translation surfaces arising from these covering constructions are known as \emph{Veech} or \emph{lattice surfaces}: they have large affine symmetry groups -- the stabilizers $SL(X, \omega)$ (known as \emph{Veech groups}) of $(X, \omega)$ under the $SL(2, \mathbb R)$-action of these surfaces are lattices. For these surfaces, Veech~\cite{veech89} showed that the existence of a saddle connection in a fixed direction (known as a \emph{cusp}) implies that the surface can be decomposed into parallel cylinders in that direction. For each cusp $c,$ we define the \emph{cusp width} relative to a group $\Gamma$ as the smallest number $w>0$ so that $\begin{pmatrix}1& w\\
0&1\end{pmatrix} \in \gamma^{-1} \Gamma \gamma$ where $\gamma(\infty)=c.$ Veech~\cite{veech} also showed that the number of saddle connections of length at most $R$ grows quadratically. We use these results to explicitly compute the asymptotic growth rate for our surfaces. We will also compute the asymptotic growth of the number of cylinders of length at most $R$, weighted by area, a quantity known as the \emph{area Siegel--Veech constant.} 

We record our main result on the Veech groups of the spectral curves of our surfaces $\Pi/\Lambda.$ Each surface results in a branched cover of the torus (branched over $0$), known as a \emph{square-tiled surface} or \emph{origami}. By \cite{gj}, these surfaces have Veech groups which are (conjugate to) finite index subgroups of the modular group $SL(2, \mathbb Z)$, the Veech group of the torus. For consistency, we apply  $\begin{pmatrix}1&\frac{1}{2}\\
0 & \frac{\sqrt{3}}{2}\end{pmatrix}^{-1}$ 
to triangle- or hexagon-tiled surfaces and consider the associated surfaces as square-tiled surfaces, and their Veech groups as subgroups of $SL(2, \mathbb Z)$. This helps to make comparisons between these surfaces easy. We abuse notation by continuing to call the surfaces by their original names. This does not change the asymptotic constants for counting, see \cite{veech}.

\begin{theorem}\label{theorem:indices}  The indices of the Veech group, the number of cusps, and the cusp widths for the translation covers of the underlying surface of Mucube, Muoctahedron, Mutetrahedron, Octa-4, and Truncated Octa-8 are given by
\begin{center}
\begin{tabular}[h]{c|c|c|c}
\hline
 & Index & Cusps & Cusp Widths \\ \hline
Mucube & $9$ & $\{\infty, 1/2, 1\}$ & $\{4, 2, 3\}$ \\ \hline
Muoctahedron & $16$ & $\{\infty, 1/5, 1/2, 1/7\}$ & $\{9, 1, 3, 3\}$ \\ \hline
Mutetrahedron & $4$ & $\{\infty, 1\}$ & $\{3, 1\}$\\ \hline 
Octa-4 & $4$ & $\{\infty, 1/2\}$ & $\{3, 1\}$ \\ \hline
Truncated Octa-8 & $15$ & $\{\infty, 3, 1, 4\}$ & $\{6, 2, 3, 4\}$ \\ \hline 
\end{tabular}

\end{center}
\end{theorem}

\noindent Given a translation surface $(X, \omega)$, we denote \begin{eqnarray}\label{eqnarray:counting}  N(R) &=& \# |\Lambda_{\omega} \cap B(0, R)|, \nonumber \\  N^{(1)}(R) &=& \# |\Lambda_{[\omega]} \cap B(0, R)|, \nonumber \\ A(R) &=& \sum_{\eta \in \La^{cyl}_{\omega} \cap B(0, R)} a(\eta), \nonumber \\ A^{(1)}(R) &=& \sum_{\eta \in \La^{cyl}_{[\omega]} \cap B(0, R)} a(\eta), \nonumber \end{eqnarray} where $[\omega]$ is the area-$1$ normalization of $\omega$ and $B(0,R)$ is the ball of radius $R$ centered at the origin. For each of our square-tiled surfaces, the area of the surface is defined as the number of squares that tile the surface. Veech~\cite{veech} showed that for the surfaces that arise in our paper, there are \emph{rational} constants $c, c^{(1)}, a, a^{(1)}$ so that 

\begin{eqnarray}\label{eqnarray:constants}  \lim_{R \rightarrow \infty} \zeta(2) \frac{N(R)}{\pi R^2} &=& c, \nonumber \\ \lim_{R \rightarrow \infty} \zeta(2) \frac{N^{(1)}(R)}{\pi R^2} &=& c^{(1)}, \nonumber \\ \lim_{R \rightarrow \infty} \zeta(2) \frac{A(R)}{\pi R^2} &=& a, \nonumber \\ \lim_{R \rightarrow \infty} \zeta(2) \frac{A^{(1)}(R)}{\pi R^2} &=& a^{(1)}.\nonumber \end{eqnarray}

\begin{theorem}\label{theorem:constants} The (normalized) asymptotic growth rates for saddle connection, cylinder, and area weighted cylinder counts for the translation covers of underlying surface of Mucube, Muoctahedron, Mutetrahedron, Octa-4, and Truncated Octa-8 are given by

\begin{center}
\begin{tabular}[h]{c|c|c|c|c}
\hline
 & $c$ & $c^{(1)}$ & $a$ & $a^{(1)}$ \\ \hline
Mucube & $1$ & $24$ & $6 \cdot 24$ & $6 \cdot 24^2$ \\ \hline
Muoctahedron & $\frac{77}{64}$ & $\frac{77}{64} \cdot 72$ & $\frac{171}{4}$ & $\frac{171}{4} \cdot 72$\\ \hline
Mutetrahedron & $\frac{19}{16}$ & $\frac{19}{16} \cdot 24$ & $15$ & $15 \cdot 24$ \\ \hline
Octa-4 & $1$ & $48$ & $6 \cdot 48$ & $6 \cdot 48^2$ \\ \hline
Truncated Octa-8 & $1$ & $120$ & $6 \cdot 120$ & $6 \cdot 120^2$ \\ \hline 
\end{tabular}
\end{center}

\end{theorem}

With the aid of the \texttt{Sage} package \cite{sage} \texttt{surface\_dynamics} \cite{dfl}, we prove that there are no closed saddle connections on any of our examples. 

\begin{theorem}\label{thm: no csc} There are no closed saddle connections on the translation covers of the underlying surface of Mucube, Muoctahedron, Mutetrahedron, Octa-4, Octa-8, or Truncated Octa-8.
\end{theorem}

Lastly, we record the constants $c, c^{(1)}, a, a^{(1)}$ for the translation covers of Platonic solids studied in Athreya--Aulicino--Hooper~\cite{AAH}.

\begin{theorem}\label{theorem:constants2} The (normalized) asymptotic growth rates for saddle connection, cylinder, and area weighted cylinder counts for the translation covers of the Platonic solids are given by:

\begin{center}
\begin{tabular}[h]{c|c|c|c|c}
\hline
 & $c$ & $c^{(1)}$ & $a$ & $a^{(1)}$ \\ \hline
Octahedron & $1$ & $12$ & $6 \cdot 12$ & $6 \cdot 12^2$ \\ \hline
Cube & $1$ & $24$ & $6 \cdot 24$ & $6 \cdot 24^2$\\ \hline
Icosahedron & $1$ & $60$ & $6 \cdot 60$ & $6 \cdot 60^2$ \\ \hline
\end{tabular}
\end{center}
\end{theorem}

The organization of this paper is as follows. In Section~\ref{sec: platonic surfaces}, we refer to \cite{thesis} and define the notion of a \emph{decoration} of a graph (Definition~\ref{defn: deco}) and show how we achieve this particular set of six polyhedral surfaces. In Section~\ref{sec: covers}, we discuss the translation covers of our examples and show that they share common covers with translation covers of Platonic solids. In Section~\ref{sec: results}, we show our results on the asymptotics of counting problems. Specifically, we will study the cusps of associated Teichm{\"u}ller curves which describe affine equivalence classes of saddle connections.  

The authors would like to thank David Aulicino, Gabriela Weitze-Schmith{\"u}sen, and the anonymous referee for the helpful comments. This material is based upon work completed while the authors were in residence at the Mathematical Sciences Research Institute in Berkeley, California, during the Fall 2019 semester. Lee was supported by the National Science Foundation under Grant No. DMS-1440140.

\section{Construction of triply periodic polyhedral surfaces} 
\label{sec: platonic surfaces}
We summarize Chapter 4 from \cite{thesis} and describe carefully the class of examples of our interest. In \cite{CP}, Coxeter and Petrie introduced three triply periodic regular polyhedral surfaces as an analogue to Platonic solids. For example, $\{4,3\}$ forms a square-tiling of the cube, $\{4,4\}$ forms a square-tiling of the plane. When one increases the valency, the faces cannot bound a convex body and their construction forces the faces to form a ``hill-valley formation'' (Figure~\ref{fig: platonic surfaces 1}). These are named Mucube, Muoctahedron, and Mutetrahedron, as they bound polyhedra that are built from multiple cubes, octahedra, and tetrahedra, respectively. In \cite{thesis}, Lee broadens this classification by allowing ``plateaus,'' while still viewing the surfaces as the boundary of a polyhedron. In other words, each surface is viewed as the boundary of a tubular neighborhood of a graph in $\R^3.$ Given a graph in $\R^3,$ Lee builds a tubular neighborhood by replacing the 0- and 1-simplices with solids and formulates a gluing pattern of the solids.

\begin{definition} \label{defn: deco} Let $\Gamma = \{V, E\}$ be a graph embedded in $\R^3$ where $V$ is a set of vertices (0-simplices) and $E$ is a set of edges (1-simplices). An edge $e \in E$ is a 2-element subset of $V$ which we denote as an unordered pair $e = \{v_1, v_2\}$ for some $v_1, v_2 \in V.$ A \textit{decoration} of $\Gamma$ is a polyhedron built by replacing the 0-simplices and 1-simplices of $\Gamma$ with convex polyhedral solids so that 1) $\Gamma$ is a deformation retract of the polyhedron and 2) the solids are identified only along faces. In essence, if a 0-simplex and a 1-simplex in $\Gamma$ are incident, then their corresponding replacement solids are identified along a common face. A \textit{regular decoration} is a decoration whose boundary surface can be denoted by Schl{\"a}fli symbols $\{p, q\}.$ An \textit{Archimedean decoration} is a decoration where the 0- and 1-simplices are replaced only by Platonic solids and Archimedean solids. 
\end{definition}

\begin{remark} We include prisms, anti-prisms, and the empty solid to replace 1-simplices but not 0-simplices. By letting empty solids replace 1-simplices, we allow two adjacent solids to retract to 0-simplices. Moreover, we will only allow the solids to be identified along one type of polygon.
\end{remark}

A graph $\Gamma$ is \textit{periodic} if $\Gamma$ is invariant under $\Lambda,$ a lattice of translations. Given a periodic graph $\Gamma,$ we define its \textit{compact quotient graph} by $\Gamma' = \{V', E'\} := \Gamma / \Lambda.$ A graph is \textit{symmetric} if given any two edges $\{v_1, v_2\}, \{v_1', v_2'\},$ there is an automorphism $\varphi: V \rightarrow V$ such that $\varphi(v_1) = v_1'$ and $\varphi(v_2) = v_2'.$ For our construction, we consider only triply periodic graphs whose quotients are symmetric. Given a quotient graph $\Gamma',$ one can define its genus by $g(\Gamma') := \textrm{e} - \textrm{v} + 1$ where $\textrm{v} = |V'|$ and $\textrm{e} = |E'|.$ 

For a given genus, there are finitely many graphs (whose quotients are symmetric) and due to the finiteness of Archimedean solids, there are only finitely many regular Archimedean decorations of that genus. We consider only those whose compact quotients are surfaces of genus-three and four that have Riemann surface structure as cyclic covers over the sphere, and are intrinsically Platonic. The surfaces are Mucube, Muoctahedron, Mutetrahedron, Octa-4, Octa-8, and Truncated Octa-8. 

The names Octa-4, Octa-8, and Truncated Octa-8 arise from the replacement solid of the 1-simplices and the valency of the graph. Adapting this notation, Mucube, Muoctahedron, and Mutetrahedron can be written as Cube-6, Truncated Octa-6, and Tetra-Truncated Tetra-4, respectively.

\section{Translation covers of Platonic surfaces}\label{sec: covers}
We devote this section to translation covers of the underlying surfaces of the six polyhedral surfaces.. 

A polyhedral surface denoted by $\{p, q\}$ is equipped with a cone metric whose cone angle is $\frac{q (p-2) \pi}{p}$ at every vertex. We take the $k$-cover so that the cover is a translation surface. 

The following table is a replication of Table 1 in \cite{AAH}.

\begin{table}[htbp]
\centering
\begin{tabular}[h]{c|c|c|c|c}
\hline
Polyhedron & Schl{\"a}fli symbols &  \begin{tabular}{@{}c@{}}Stratum of the \\ $k$-differential\end{tabular} & \begin{tabular}{@{}c@{}} Stratum of the \\ translation cover\end{tabular} & \begin{tabular}{@{}c@{}}Genus of the \\ translation cover\end{tabular} \\ \hline
Mucube & $\{4, 6\}$ & $\mathcal{H}_2(1^8)$ & $\mathcal{H}_1(2^8)$ & 9 \\ \hline
Muoctahedron & $\{6, 4\}$ & $\mathcal{H}_3(1^{12})$ & $\mathcal{H}_1(3^{12})$ & 19 \\ \hline
Mutetrahedon & $\{6, 6\}$ & $\mathcal{H}_2(1^4)$ & $\mathcal{H}_1(1^8)$ & 5 \\ \hline 
Octa-4 & $\{3, 8\}$ & $\mathcal{H}_3(1^{12})$ & $\mathcal{H}_1(3^{12})$ & 19 \\ \hline
Octa-8 & $\{3, 12\}$ & $\mathcal{H}_1(1^6)$ & $\mathcal{H}_1(1^6)$ & 4 \\ \hline
Truncated Octa-8 & $\{4, 5\}$ & $\mathcal{H}_4(1^{24})$ & $\mathcal{H}_1(4^{24})$ & 49 \\ \hline
\end{tabular}
	\label{tab: platonic surfaces}
	\caption{Strata of polyhedral surfaces and their translation covers.}
\end{table}

In \cite{thesis}, Lee finds various holomorphic maps from the quotient of Octa-8 to the sphere to obtain a basis of holomorphic 1-forms. In the following theorem, we show that the associated 1-form corresponds to the unfolding of the octahedron. 

\begin{theorem}
The underlying Riemann surface structure on Octa-8 corresponds to the translation cover of the octahedron.\end{theorem}

\begin{proof} The underlying surface of Octa-8 is a genus-four compact Riemann surface with identification of edges as described in Figure~\ref{fig: octa8_hyp}. 

\begin{figure}[htbp] 
	\centering
	\includegraphics[width=2.5in]{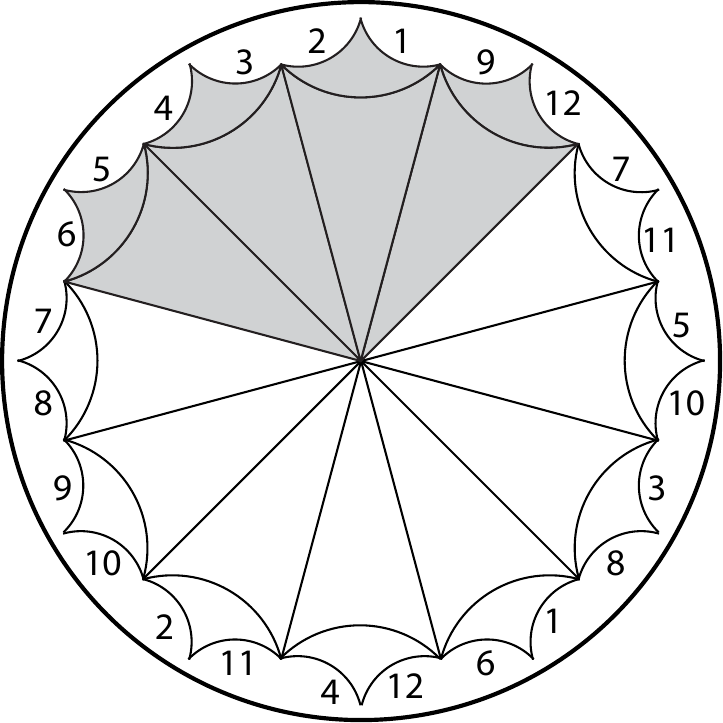}
	\caption{Hyperbolic description of the fundamental piece of Octa-8. Adapted from \cite{thesis}}	
	\label{fig: octa8_hyp} 
\end{figure}

The surface is invariant under the order-twelve rotation $\varphi$ about the center of the tessellation. The quotient $X / \langle \varphi \rangle$ is a sphere where the covering is branched over three points. However, $X / \langle \varphi^4 \rangle$ is also a sphere (shaded region in Figure~\ref{fig: octa8_hyp}) and the threefold map is branched over six points. In Section 5.2 of \cite{thesis}, Lee shows that the quotient sphere $X / \langle \varphi^4 \rangle$ is conformally equivalent to a regular octahedron. In other words, the exterior derivative of this particular map is a holomorphic 1-form corresponding to the translation cover of the octahedron.

\end{proof}

Next we show that the cube (quartic differential) and the compact quotient of Mucube (quadratic differential) share the same translation cover. 

\begin{theorem} 
The polyhedral metric on the underlying surface of Mucube yields a quadratic differential that corresponds to the half-translation cover of the cube.\end{theorem}

\begin{proof} We prove this theorem by picture. Figure~\ref{fig: two_cube} describes Mucube and the twofold cover of the cube where the edges are identified by translation or a half-translation. They are identical as half-translation surfaces, hence their twofold translation covers are identical.

\begin{figure}[htbp] 
	\centering
    \includegraphics[width=4.5in]{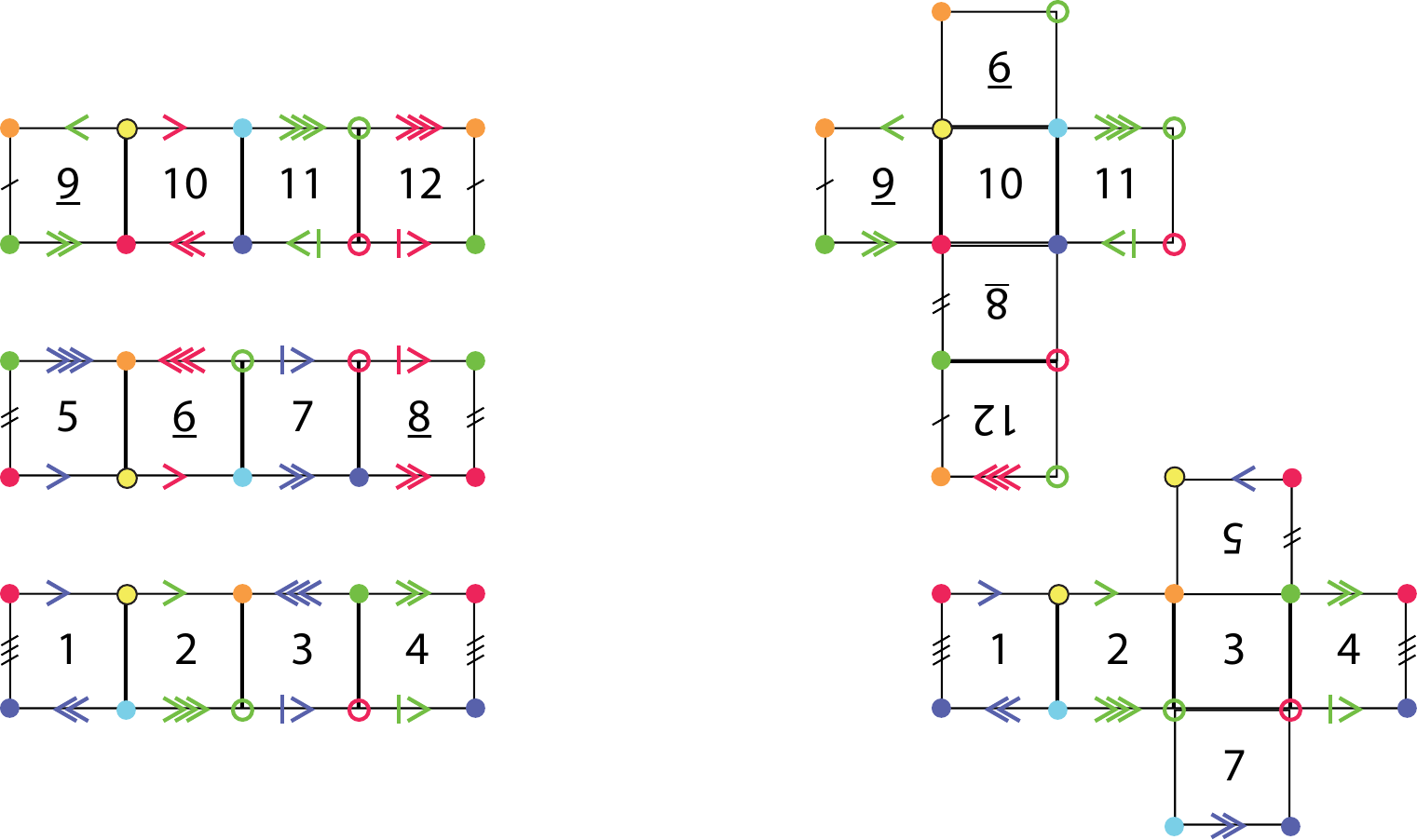}
    \caption{Mucube (left) and the twofold cover of the cube (right).}	
    \label{fig: two_cube}
    \end{figure}
\end{proof}

\section{Asymptotics of counting problems}
\label{sec: results}

\begin{figure}[htbp] 
	\centering
	\includegraphics[width=2.5in]{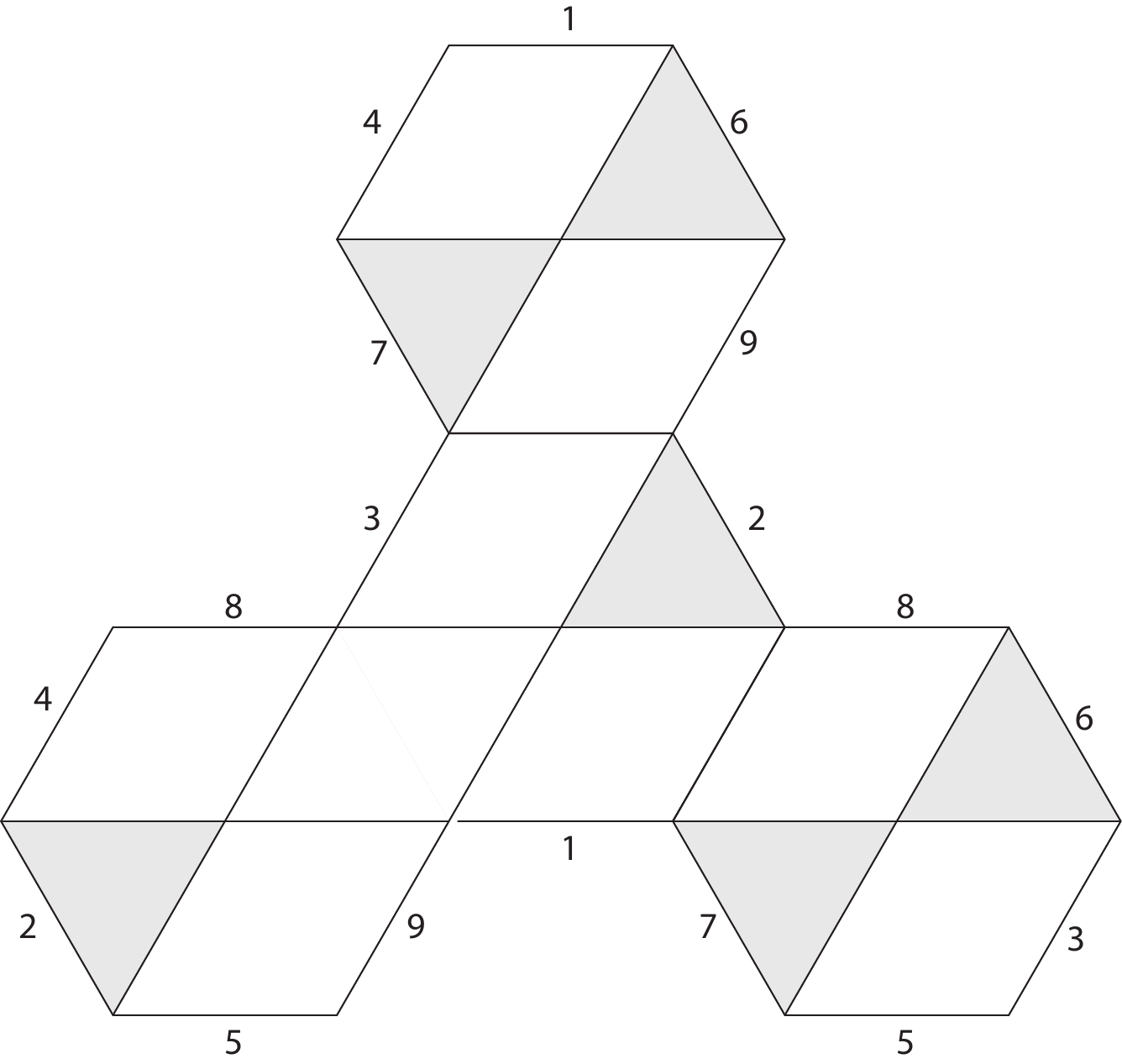}
	\caption{Subdivision of hexagons on Mutetrahedron}	
	\label{fig: mutetra_subdivision}
\end{figure}

In the following subsections we study the Teichm{\"u}ller curve associated to the translation cover of each polyhedral surface in Table~\ref{tab: platonic surfaces}. To ease our computations, we describe each translation surface as a surface tiled by unit squares. The advantage of this tool is that a surface is completely defined by a horizontal permutation $\sigma_h$ and a vertical permutation $\sigma_v$ on the squares. The top of square $i$ is glued to the bottom of square $\sigma_v(i)$, and the right side of square $i$ is glued to the left side of square $\sigma_h(i)$. The polyhedral surfaces tiled by triangles cover the doubled triangle (a rhombus), which we map to a square as shown in \cite{AAH}, and as discussed earlier. Our examples include surfaces tiled by hexagons (Muoctahedron and Mutetrahedron), which can also be described as square-tiled surfaces. We subdivide each hexagon into six triangles and then pair two triangles to form a rhombus. Figure~\ref{fig: mutetra_subdivision} illustrates the subdivision of Mutetrahedron into rhombi where the numbered edges are identified by half-translations. Note that the covering is not regular over rhombi.

The input to the \texttt{Sage} package \texttt{surface\_dynamics} is the two permutations that define a square-tiled surface. It yields Theorem~\ref{theorem:indices} where it computes the Veech group of the Teichm{\"u}ller curve and its cusps and cups widths. We will use Theorem~\ref{theorem:indices} to study the Teichm{\"u}ller curves associated to the polyhedral surfaces. 

\subsection*{Mutetrahedron}
The underlying surface of Mutetrahedron is tiled by four hexagons. The cone angles at each vertex is $2 \pi,$ however the edges are identified by translation and a $180^{\circ}$-rotation. Hence its spectral curve is a 2-cover. By subdivision of hexagons into squares, the translation cover of Mutetrahedron is tiled by 24 squares, described by a horizontal permutation $\sigma_h$ and a vertical permutation $\sigma_v.$ \begin{align}\sigma_h = (1,2,3,4,5,6)(7,8,9,10,11,12)(13,14,15,16,17,18)(19,20,21,22,23,24)\\
\sigma_v = (1,7,13,11,3,21)(2,20,14,12,18,22)(4,10,16,8,6,24)(5,23,17,9,15,19)\end{align} 

\begin{figure}[htbp] 
	\centering
	\includegraphics[width=5in]{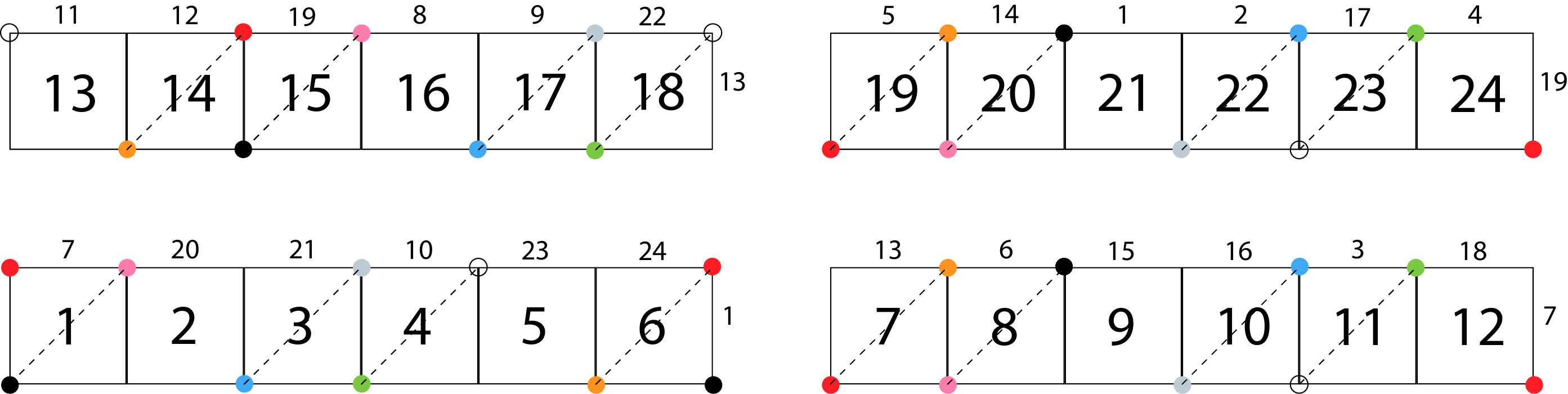}
	\caption{Double cover of the underlying surface of Mutetrahedron.}	
	\label{fig: csc_on_MT1}
\end{figure}

The translation cover is a genus five surface with eight simple zeros. Its Veech group $\Gamma$ is an index-4 subgroup of $SL(2,\Z),$ which is described completely by the action of the generators of $SL(2,\Z)$ on the right cosets $\Gamma \backslash SL(2,\Z).$ Given the generators $s_2 = \bigl( \begin{smallmatrix}0 & -1\\1 & 0\end{smallmatrix}\bigr),$ $s_3 = \bigl( \begin{smallmatrix}0 & 1\\-1 & 1\end{smallmatrix}\bigr),$ $l = \bigl( \begin{smallmatrix}1 & 1\\0 & 1\end{smallmatrix} \bigr),$ and $r = \bigl( \begin{smallmatrix}1 & 0\\1 & 1\end{smallmatrix} \bigr)$ of $SL(2,\Z),$ $\Gamma$ has the following representation:

$$\begin{array}{l}
s_2=(1,4)(2,3),\\
s_3=(1,3,2),\\
l=(1,3,4),\\
r=(1,2,4),\end{array}$$

where the cosets are identified with elements in $\{1,2,3,4\}.$ In this and future examples, 1 corresponds to the coset containing the identity.

The cusps are at $\infty$ and 1 with width 3 and 1, respectively. Each cusp corresponds to a rational direction in which the surface decomposes into parallel cylinders. We denote the saddle connections by $\mathbf{v}_{\infty} = \bigl[ \begin{smallmatrix}0\\1\end{smallmatrix} \bigr],$ $\mathbf{v}_{\infty}' = \bigl[ \begin{smallmatrix}0\\2\end{smallmatrix} \bigr],$ and 
$\mathbf{v}_1 = \bigl[ \begin{smallmatrix}1\\1\end{smallmatrix} \bigr].$ However, neither of these are closed saddle connections as they do not connect vertices of the same color (Figure~\ref{fig: csc_on_MT1}). To compute the quadratic asymptotics of saddle connections, we follow Section 16 of Veech~\cite{veech} and Appendix C of Athreya--Chaika--Leli\`{e}vre~\cite{ACL}. Veech~\cite[Theorem 16.1]{veech} showed that for any non-uniform lattice $\Gamma \subset SL(2, \R)$ and any vector $\mathbf{v}$ stabilized by a maximal parabolic subgroup $\Lambda \subset \Gamma,$ we have $$\lim_{R \rightarrow \infty} \frac{|g \Gamma \mathbf{v} \cap B(0,R)|}{\pi R^2} = c(\Gamma, \mathbf{v}).$$ Our goal is to find this limit for all saddle connection vectors $\mathbf{v} \in \Lambda_{\omega},$ that is, $$\lim\limits_{R \rightarrow \infty} \frac{N(R)}{\pi R^2} = \sum\limits_{\mathbf{v} \in \Lambda_{\omega}} c(\Gamma, \mathbf{v}).$$ Given $\mathbf{v},$ we compute $c(\Gamma, \mathbf{v})$ as follows. Let $g_0 \in SL(2,\R)$ be such that $g_0^{-1} \Lambda g_0 = \Lambda_0,$ with $\Lambda_0 = \left\{\bigl(\begin{smallmatrix}1 & n\\
0 & 1\end{smallmatrix}\bigr) : n \in \Z \right\}.$ Setting $\mathbf{v}_0 = 1$ and  $g_0 \mathbf{v}_0 = t \mathbf{v},$ we have $c(\Gamma, \mathbf{v}) = t^2 \textrm{vol}(\mathbb{H}^2/\Gamma)^{-1}.$ Here, $t^2$ corresponds to the width of the cusp. For $\mathbf{v} = \mathbf{v}_{\infty},$ we have $$g_0 = \begin{pmatrix}\sqrt{3} & 0 \\
0 & 1/\sqrt{3}\end{pmatrix}$$ hence $c(\Gamma, \mathbf{v})= 3 \cdot \frac{6}{\pi^2} \cdot \frac{1}{4} = \frac{9}{2 \pi^2}.$ 
Since $\mathbf{v}_{\infty}' = 2 \mathbf{v}_{\infty},$ we have $c(\Gamma, \mathbf{v}_{\infty}') = \left(\frac{1}{2}\right)^2 c(\Gamma, \mathbf{v}_{\infty}).$ 

\noindent Following this computation, we have $$\lim_{R \rightarrow \infty} \dfrac{N(R)}{\pi R^2} = \left(\frac{3}{4} \left(1 + \frac{1}{4}\right) + \frac{1}{4}\right) \frac{6}{\pi^2} = \frac{19}{16} \cdot \frac{6}{\pi^2}.$$

\noindent We normalize the area of the surface to $1$ and denote by $N^{(1)}(R)$ the number of saddle connection vectors of length at most $R$ on the unit-area surface. Then $$\lim_{R \rightarrow \infty} \dfrac{N^{(1)}(R)}{\pi R^2} = \frac{19}{16} \cdot \frac{6}{\pi^2} \cdot 24.$$ Now we consider the cylinder decomposition of the surface arising from each saddle connection vector and compute the asymptotic growth of (weighted) cylinders whose core curves are at most length $R.$ In the vertical direction, the surface decomposes into sixteen cylinders: eight have $\mathbf{v}_{\infty}$ as their core curves and eight have $\mathbf{v}_{\infty}'$ as their core curves, hence the area of the cylinders are 1 and 2 respectively. Similarly, in the direction of $\mathbf{v}_1,$ the surface decomposes into eight cylinders of unit area and eight cylinders of area 2. Hence $$\lim\limits_{R \rightarrow \infty} \dfrac{A(R)}{\pi R^2} = \left(\dfrac{3}{4} \left(8 + \dfrac{16}{4}\right) + \dfrac{1}{4}(8 + 16)\right) \dfrac{6}{\pi^2} = 15 \cdot \dfrac{6}{\pi^2}$$ and on the unit-area surface, we have $\lim\limits_{R \rightarrow \infty} \dfrac{A^{(1)}(R)}{\pi R^2} = 15 \cdot \dfrac{6}{\pi^2} \cdot 24.$

\subsection*{Mucube} 
The translation cover of the underlying surface of Mucube is defined by the following horizontal and vertical permutations on 24 squares. It is a genus nine surface with eight simple zeros.
\begin{align}\sigma_h = (1,2,3,4)(5,6,7,8)(9,10,11,12)(13,14,15,16)(17,18,19,20)(21,22,23,24)\\
\sigma_v = (1,9,14,22)(2,20,13,7)(3,24,16,11)(4,5,15,18)(6,10,17,21)(8,23,19,12)\end{align} 

\begin{figure}[htbp] 
	\centering
	\includegraphics[width=5in]{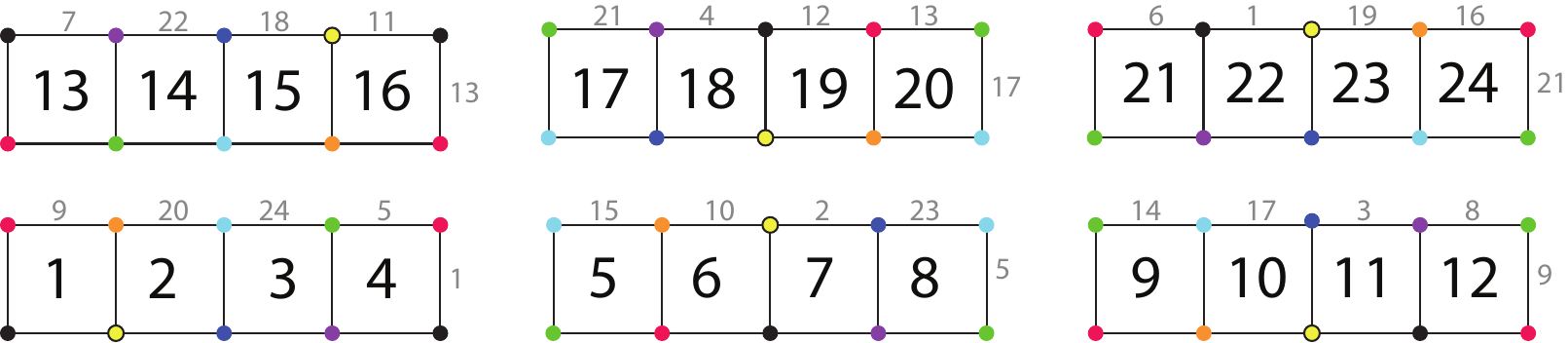}
	\caption{Double cover of the underlying surface of Mucube}	
	\label{fig: 2cover_of_prim}
\end{figure}

The Veech group is an index-9 subgroup of $SL(2,\Z)$ which we describe by the action of the generators of $SL(2,\Z)$ on the right cosets: $$\begin{array}{l}s_2=(2,4)(3,7)(5,9)(6,8)\\
s_3=(1,9,6)(2,5,7)(3,8,4)\\
l=(1,8,7,9)(2,3)(4,6,5)\\
r=(1,5,3,6)(2,9,8)(4,7)\end{array}$$ 

The cusps are at $\infty,$ 1/2, and 1, with cusp width 4, 2, and 3, respectively. These are associated to the vertical direction, direction of slope $1/2,$ and $1/3.$ Note that none of the saddle connections are closed as there is no vertex on both top and bottom of any cylinder in the cusp directions (Figure~\ref{fig: 2cover_of_prim}). Then we have the following quantities: $$\lim_{R \rightarrow \infty} \dfrac{N(R)}{\pi R^2} = \left(\frac{4}{9} + \frac{2}{9} + \frac{3}{9}\right) \frac{6}{\pi^2} = \frac{6}{\pi^2} \qquad \textrm{and} \qquad \lim_{R \rightarrow \infty} \dfrac{N^{(1)}(R)}{\pi R^2} = \frac{6 \cdot 24}{\pi^2}.$$

Since every cylinder has unit-area, 
$$\lim_{R \rightarrow \infty} \frac{A(R)}{\pi R^2} = \frac{6 \cdot 24}{\pi^2} \qquad \textrm{and} \qquad \lim_{R \rightarrow \infty} \frac{A^{(1)}(R)}{\pi R^2} = \frac{6 \cdot 24^2}{\pi^2}.$$

\subsection*{Muoctahedron} 
The translation cover of the underlying surface of Muoctahedron is a genus-nineteen surface with twelve order-three zeros. The following vertical and horizontal permutations define the square-tiled surface.

$$\begin{array}{ll}\sigma_v = & (1,29,20,36,39,54,37,3,10)(2,18,46,38,53,21,28,19,11)(4,60,70,58,35,23,33,6,16)\\
& (5,15,24,34,22,71,59,69,17)(7,42,51,40,57,64,55,9,13)(8,14,50,41,52,72,56,65,12)\\
& (25,31,27,66,62,68,47,44,49)(26,32,43,48,45,61,67,63,30)]\\
\sigma_h = & (1,2,3,4,5,6,7,8,9)(10,11,12,13,14,15,16,17,18)(19,20,21,22,23,24,25,26,27)\\
& (28,29,30,31,32,33,34,35,36)(37,38,39,40,41,42,43,44,45)(46,47,48,49,50,51,52,53,54)\\
& (55,56,57,58,59,60,61,62,63)(64,65,66,67,68,69,70,71,72)\end{array}$$ 

\begin{figure}[htbp] 
	\centering
	\includegraphics[width=5in]{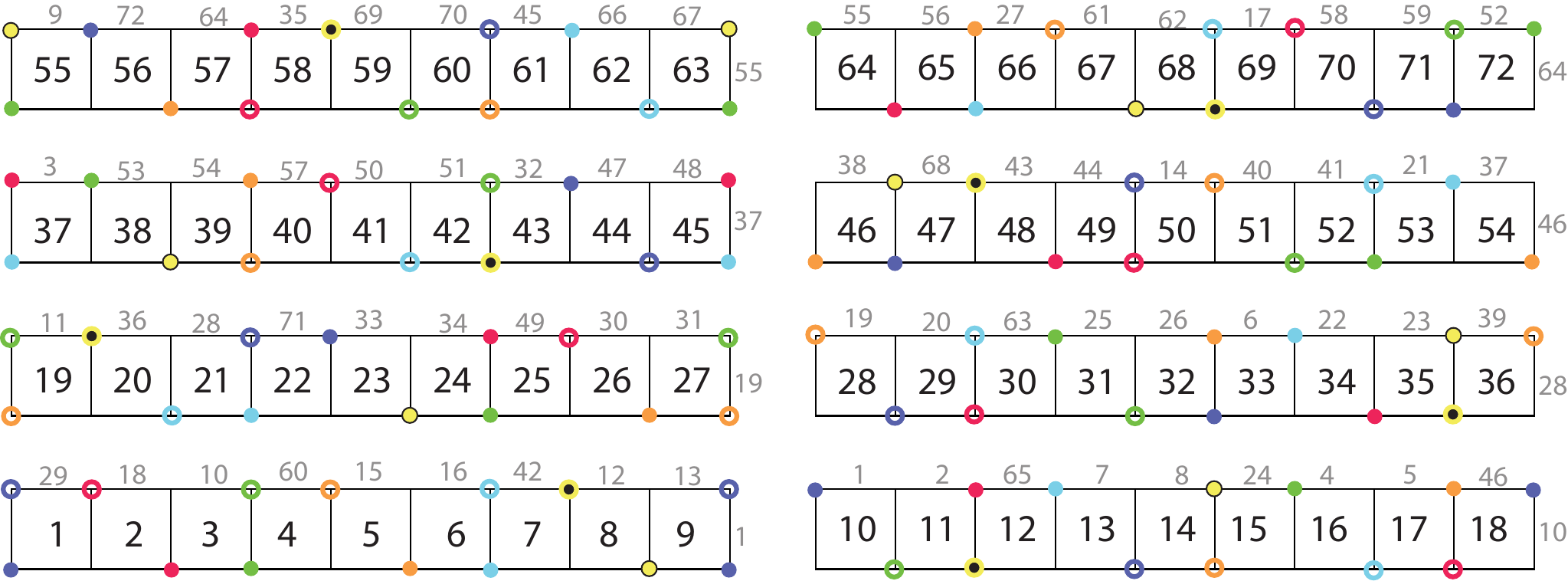}
	\caption{Triple cover of the underlying surface of Muoctahedron}	
	\label{fig: 3cover_of_octaprim_in square form}
\end{figure}

The Veech group is an index-16 subgroup of $SL(2,\Z),$ which we describe by the action of generators on the right cosets:
$$\begin{array}{l}s_2=(1,14)(2,9)(3,11)(4,10)(5,6)(7,15)(8,12)(13,16)\\
s_3=(2,16,9)(3,15,6)(4,14,11)(5,8,10)(7,13,12)\\
l=(1,14,10,12,16,9,13,15,11)(3,5,4)(6,7,8)\\
r=(1,14,3,7,16,2,13,8,4)(5,12,15)(6,11,10)\end{array}$$

The cusps are at $\infty,$ 1/5, 1/2, and 1/7, with width 9, 1, 3, and 3, respectively. Again, none of the saddle connections are closed. The number of saddle connection vectors in the disk of radius $R$ grows quadratically and we get $$\lim_{R \rightarrow \infty} \dfrac{N(R)}{\pi R^2} = \left(\frac{9}{16} \left(1 + \frac{1}{4}\right) + \frac{1}{16} \left(1 + \frac{1}{4}\right) + \frac{3}{16} \left(1 + \frac{1}{4}\right) + \frac{3}{16}\right) \frac{6}{\pi^2}  = \frac{77}{64} \cdot \frac{6}{\pi^2}$$ and $$\lim_{R \rightarrow \infty} \dfrac{N^{(1)}(R)}{\pi R^2} = \frac{77}{64} \cdot \frac{6}{\pi^2} \cdot 72.$$

For all directions, there are 24 cylinders of unit-area and 24 cylinders of area 2. For directions $\infty,$ 1/5, and 1/2, the area of the cylinder is proportional to the length of the core curve. 

$$\begin{array}{ll}\lim\limits_{R \rightarrow \infty} \dfrac{A(R)}{\pi R^2} & = \left(\dfrac{9}{16} \left(24 + \dfrac{48}{4}\right) + \dfrac{1}{16} \left(24 + \dfrac{48}{4}\right) + \dfrac{3}{16} \left(24 + \dfrac{48}{4}\right) + \dfrac{3}{16} \cdot (24 + 48)\right) \dfrac{6}{\pi^2}\\
&\\
& = \dfrac{171}{4} \cdot \dfrac{6}{\pi^2}\end{array}$$ and 
$$\lim_{R \rightarrow \infty} \frac{A^{(1)}(R)}{\pi R^2} = \frac{171}{4} \cdot \frac{6}{\pi^2} \cdot 72.$$

\subsection*{Octa-4} 
\begin{figure}[htbp] 
	\centering
	\includegraphics[width=5in]{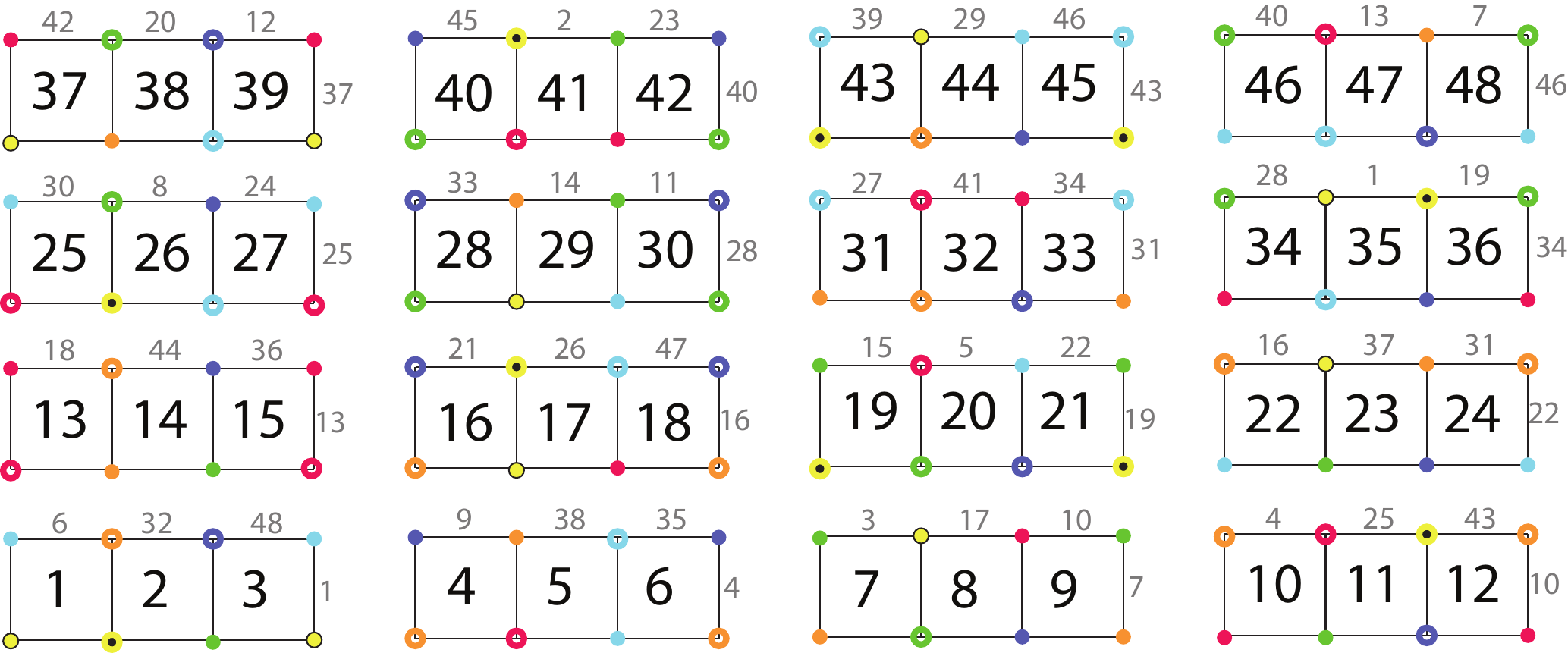}
	\caption{Triple cover of the underlying surface of Octa-4}	
	\label{fig: 3cover_of_octa4}
\end{figure}

The translation cover of the underlying surface of Octa-4 is a genus-nineteen surface with twelve order-three zeros. The cover is defined by vertical and horizontal permutations as described in Figure~\ref{fig: 3cover_of_octa4}.

$$\begin{array}{ll}\sigma_v = & (1,6,35)(2,32,41)(3,48,7)(4,9,10)(5,38,20)(8,17,26)(11,25,30)(12,43,39)\\
& (13,18,47)(14,44,29)(15,36,19)(16,21,22)(23,37,42)(24,31,27)(28,33,34)(40,45,46) \\
\sigma_h = & (1,2,3)(4,5,6)(7,8,9)(10,11,12)(13,14,15)(16,17,18)(19,20,21)(22,23,24)\\
& (25,26,27)(28,29,30)(31,32,33)(34,35,36)(37,38,39)(40,41,42)(43,44,45)(46,47,48)\end{array}$$

Its Veech group is an index-four subgroup of $SL(2,\Z)$ with the following description on the right cosets: $$s_2=(1,2)(3,4), \qquad s_3=(1,3,4), \qquad l=(1,3,2), \qquad r=(1,4,2).$$

The two cusps are $\infty$ and $1/2$ with width 3 and 1, respectively. There are no closed saddle connections. We achieve $$\lim_{R \rightarrow \infty} \dfrac{N(R)}{R^2} = \left(\frac{3}{4} + \frac{1}{4}\right) \frac{6}{\pi^2} = \frac{6}{\pi^2} \qquad \textrm{and} \qquad \lim_{R \rightarrow \infty} \dfrac{N^{(1)}(R)}{\pi R^2} = \frac{6 \cdot 48}{\pi^2}.$$

Since all cylinders are unit-area, we get $$\lim_{R \rightarrow \infty} \frac{A(R)}{\pi R^2} = \frac{6 \cdot 48}{\pi^2} \qquad \textrm{and} \qquad \lim_{R \rightarrow \infty} \frac{A^{(1)}(R)}{\pi R^2} = \frac{6 \cdot 48^2}{\pi^2}.$$

\subsection*{Truncated Octa-8}
The translation cover of the underlying surface of Truncated Octa-8 is a genus-49 surface with 24 order-4 zeros. The square-tiled surface is defined by the following permutations on 120 squares. 

$$\begin{array}{ll}\sigma_v = & (1,34,75,81,120,7)(2,20,111,80,62,25)(3,40,67,79,107,15)(4,10,117,84,78,31)\\
& (5,28,65,83,114,23)(6,18,104,82,70,37)(8,52,71,74,96,17)(9,24,86,73,64,43)\\
& (11,14,93,77,68,49)(12,46,61,76,89,21)(13,22,99,69,66,55)(16,58,63,72,102,19)\\
& (26,45,115,110,90,33)(27,38,98,109,103,56)(29,36,87,113,118,48)(30,59,106,112,101,41)\\
& (32,92,108,116,50,39)(35,42,53,119,105,95)(44,57,91,85,97,51)(47,54,100,88,94,60)\\
\sigma_h = & (1,2,3,4,5,6)(7,8,9,10,11,12)(13,14,15,16,17,18)(19,20,21,22,23,24)\\
& (25,26,27,28,29,30)(31,32,33,34,35,36)(37,38,39,40,41,42)(43,44,45,46,47,48)\\
& (49,50,51,52,53,54)(55,56,57,58,59,60)(61,62,63,64,65,66)(67,68,69,70,71,72)\\
& (73,74,75,76,77,78)(79,80,81,82,83,84)(85,86,87,88,89,90)(91,92,93,94,95,96)\\
& (97,98,99,100,101,102)(103,104,105,106,107,108)(109,110,111,112,113,114)\\
& (115,116,117,118,119,120)\end{array}$$

\begin{figure}[htbp] 
	\centering
	\includegraphics[width=5in]{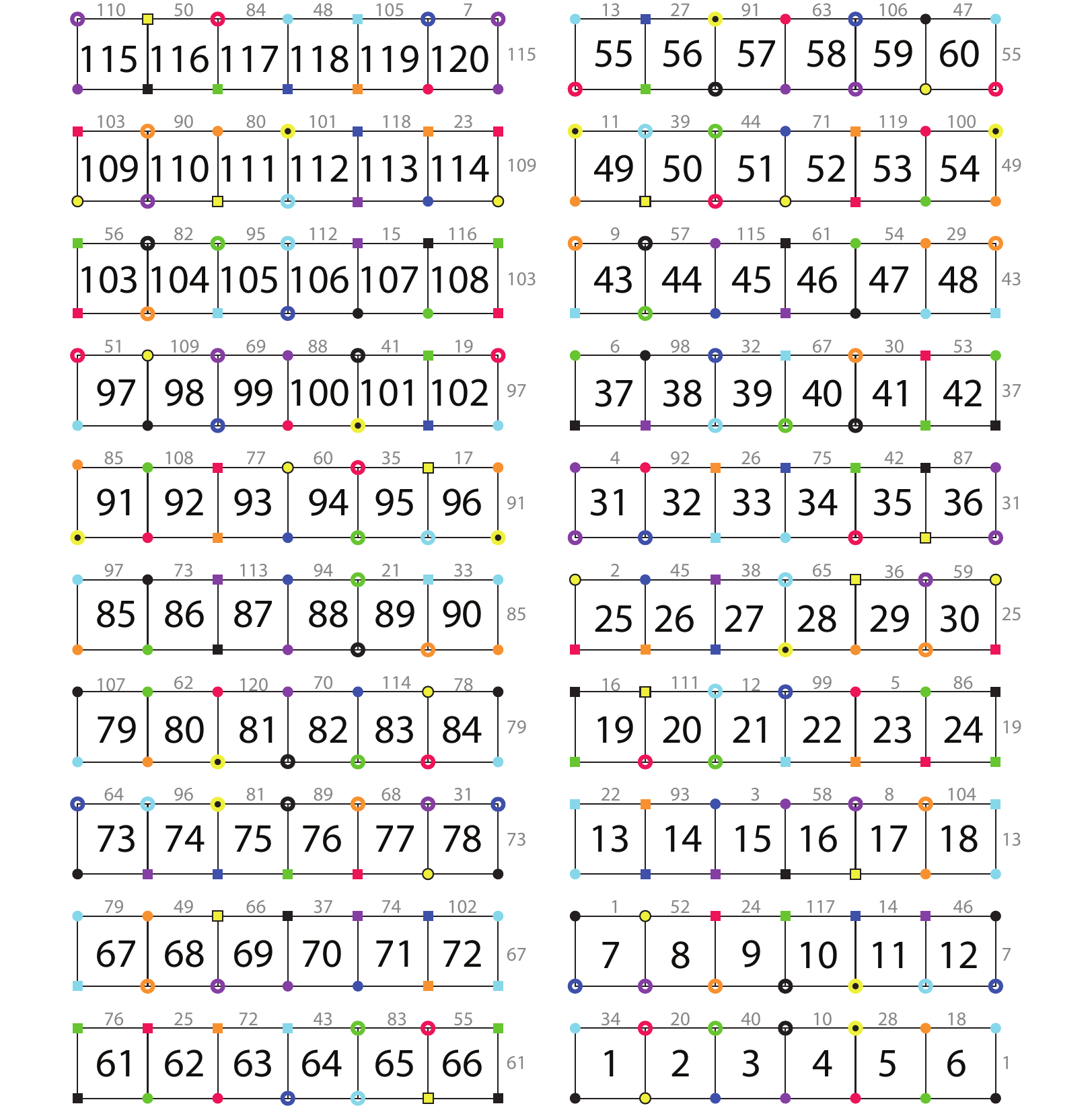}
	\caption{Quadruple cover of the underlying surface of Truncated Octa-8}	
	\label{fig: 4cover_of_bring}
\end{figure}

Its Veech group is an index-15 subgroup of $SL(2,\Z)$ with the following description on the right cosets:
$$\begin{array}{l}s_2=(2,7)(3,11)(4,8)(5,13)(6,12)(9,15)(10,14)\\
s_3=(1,13,6)(2,9,11)(3,14,7)(4,10,12)(5,15,8)\\
l=(1,12,14,11,15,13)(2,3)(4,6,5)(7,10,8,9)\\
r=(1,5,9,3,10,6)(2,15,4,14)(7,11)(8,13,12)\end{array}.$$

The cusps are at $\infty,$ 1/3, 1, and 1/4, and none of the saddle connections are closed. We get $$\lim_{R \rightarrow \infty} \dfrac{N(R)}{\pi R^2} = \left( \frac{6}{15} + \frac{2}{15} + \frac{3}{15} + \frac{4}{15}\right) \frac{6}{\pi^2} = \frac{6}{\pi^2} \qquad \textrm{and} \qquad \lim_{R \rightarrow \infty} \dfrac{N^{(1)}(R)}{\pi R^2} = \frac{6 \cdot 120}{\pi^2}.$$

Since every cylinder is unit-area, 
$$\lim_{R \rightarrow \infty} \frac{A(R)}{\pi R^2} = \frac{6 \cdot 120}{\pi^2} \qquad \textrm{and} \qquad \lim_{R \rightarrow \infty} \frac{A^{(1)}(R)}{\pi R^2} = \frac{6 \cdot 120^2}{\pi^2}.$$
 
\subsection*{Platonic solids}
In this section, we present the quadratic asymptotics of saddle connection vectors on the translation cover of Platonic solids represented as square-tiled surfaces. The following is Table 2 from \cite{AAH}. 

\begin{table}[htbp]
\centering
\begin{tabular}[h]{c|c|c}
\hline
 & Cusps & Cusp Widths \\ \hline
Tetrahedron & $\{\infty\}$ & $\{1\}$ \\ \hline
Octahedron & $\{\infty, 1\}$ & $\{3, 1\}$ \\ \hline

Cube & $\{\infty, 1/2, 1\}$ & $\{4, 2, 3\}$ \\ \hline

Icosahedron & $\{\infty, 1/3, 1/2\}$ & $\{5, 2, 3\}$ \\ \hline

\end{tabular}
	\label{tab: cusp and cusp widths}
	\caption{Cusp and cusp widths of Platonic solids}
\end{table}

Then we have the following quantities for the spectral curve of the octahedron, cube, and icosahedron:

\subsubsection*{3-cover of the octahedron}
$$\lim_{R \rightarrow \infty} \dfrac{N(R)}{\pi R^2} = \left(\frac{3}{4} + \frac{1}{4}\right) \frac{6}{\pi^2} = \frac{6}{\pi^2} \qquad \textrm{and} \qquad \lim_{R \rightarrow \infty} \dfrac{N^{(1)}(R)}{\pi R^2} = \frac{6}{\pi^2} \cdot 12,$$

and all cylinders are unit-area: $$\lim_{R \rightarrow \infty} \frac{A(R)}{\pi R^2} = \frac{6}{\pi^2} \cdot 12 \qquad \textrm{and} \qquad \lim_{R \rightarrow \infty} \frac{A^{(1)}(R)}{\pi R^2} = \frac{6}{\pi^2} \cdot 12^2.$$

\subsubsection*{4-cover of the cube (identical to the 2-cover of Mucube)}
$$\lim_{R \rightarrow \infty} \dfrac{N(R)}{\pi R^2} = \left(\frac{4}{9} + \frac{2}{9} + \frac{3}{9}\right) \frac{6}{\pi^2} = \frac{6}{\pi^2} \qquad \textrm{and} \qquad \lim_{R \rightarrow \infty} \dfrac{N^{(1)}(R)}{\pi R^2} = \frac{6 \cdot 24}{\pi^2}.$$

Similarly, all cylinders are unit-area, hence $$\lim_{R \rightarrow \infty} \frac{A(R)}{\pi R^2} = \frac{6 \cdot 24}{\pi^2} \qquad \textrm{and} \qquad \lim_{R \rightarrow \infty} \frac{A^{(1)}(R)}{\pi R^2} = \frac{6 \cdot 24^2}{\pi^2}.$$

\subsubsection*{6-cover of the icosahedron}

$$\lim_{R \rightarrow \infty} \dfrac{N(R)}{R^2} = \left(\frac{5}{10} + \frac{2}{10} + \frac{3}{10}\right) \frac{6}{\pi^2} = \frac{6}{\pi^2} \qquad \textrm{and} \qquad \lim_{R \rightarrow \infty} \dfrac{N^{(1)}(R)}{\pi R^2} = \frac{6 \cdot 60}{\pi^2}.$$
Also, all cylinders are unit-area, hence 
$$\lim_{R \rightarrow \infty} \frac{A(R)}{\pi R^2} = \frac{6 \cdot 60}{\pi^2} \qquad \textrm{and} \qquad  \lim_{R \rightarrow \infty} \frac{A^{(1)}(R)}{\pi R^2} = \frac{6 \cdot 60^2}{\pi^2}.$$



\end{document}